\documentclass[12pt,twoside]{amsart}
\usepackage{amssymb}
\usepackage[all]{xy}

\nonstopmode

\numberwithin{equation}{section} 

\font\tengothic=eufm10 scaled\magstep 1 \font\sevengothic=eufm7
scaled\magstep 1
\newfam\gothicfam
     \textfont\gothicfam=\tengothic
     \scriptfont\gothicfam=\sevengothic

\newtheorem{theorem}{Theorem}[section]
\newtheorem{lemma}[theorem]{Lemma}
\newtheorem{proposition}[theorem]{Proposition}
\newtheorem{corollary}[theorem]{Corollary}

\theoremstyle{definition}
\newtheorem{definition}[theorem]{Definition} 
\newtheorem{remark}[theorem]{Remark}
\newtheorem{problem}[theorem]{Problem}
\newtheorem{example}[theorem]{Example}

\newtheorem{notation}[theorem]{Notation}

\newcommand{\Hom}{\operatorname{Hom}}
\newcommand{\Ext}{\operatorname{Ext}}
\newcommand{\Supp}{\operatorname{Supp}}
\newcommand{\Pic}{\operatorname{Pic}}
\newcommand{\rank}{\operatorname{rank}}

\newcommand{\cE}{{\mathcal E}}

\newcommand{\cO}{{\mathcal O}}
\newcommand{\cL}{{\mathcal L}}

\newcommand{\cP}{{\mathcal P}}

\newcommand{\cD}{{\mathcal D}}

\newcommand {\RR}{\mathbb{R}}

\newcommand {\ZZ}{\mathbb{Z}}

\newcommand {\PP}{\mathbb{P}}

\begin{document}
\title[Derived category of toric fibrations]
{Derived category of toric fibrations}

\author[L.\ Costa, S. \ Di Rocco, R.M.\ Mir\'o-Roig]{L.\ Costa$^*$, S. \ Di Rocco$^{**}$, R.M.\
Mir\'o-Roig$^{*}$}

\address{Facultat de Matem\`atiques,
Departament d'Algebra i Geometria, Gran Via de les Corts Catalanes
585, 08007 Barcelona, SPAIN } \email{costa@ub.edu}
\address{Department of Mathematics, KTH, SE-10044 Stockholm, Sweden }
\email{dirocco@math.kth.se}
\address{Facultat de Matem\`atiques,
Departament d'Algebra i Geometria, Gran Via de les Corts Catalanes
585, 08007 Barcelona, SPAIN } \email{miro@ub.edu}

\date{\today}
\thanks{$^*$ Partially supported by MTM2010-15256.}
\thanks{$^{**}$ Partially supported by  Vetenskapsr{\aa}det's grant
NT:2006-3539.}



\begin{abstract} In this paper  we give a structure theorem for the derived category $D^b(X)$
 of a toric fiber bundle $X$ over $Z$ with fiber $F$ provided that $F$ and $Z$ have both a full strongly
  exceptional collection of line bundles.
\end{abstract}



\maketitle



\section{Introduction} \label{intro}
Let $X$ be a smooth projective variety defined over an
algebraically closed field $K$ of characteristic zero and let
$D^b(X)$ be the derived category of bounded complexes of coherent
sheaves of ${\cO}_X$-modules. $D^b(X)$ is one of the most
important algebraic invariants of a smooth projective variety $X$
and, in spite of the increasing interest in understanding the
structure of $D^b(X)$, very little progress has been achieved. The
study of $D^b(X)$ dates back to the late 70's when Beilinson
described the derived category of projective spaces (\cite{Be})
and it has became one of the most important topics in Algebraic
Geometry. Beside other reasons this is due to the Homological
Mirror Symmetry Conjecture of Kontsevich \cite{Ko} which states
that there is an equivalence of categories between the derived
category of coherent sheaves on a Calabi-Yau  variety and the
derived Fukai category of its mirror.

\vskip 2mm

An important approach to determine the structure of $D^b(X)$
 is to construct tilting bundles.
Following terminology of representation
theory (cf. \cite{Ba}) a coherent sheaf $T$ of ${\cO}_X$-modules on a
smooth projective variety $X$ is called a {\em tilting sheaf} (or,
when it is locally free, a {\em tilting bundle}) if
\begin{itemize}
\item[(i)] it has no higher self-extensions, i.e.
$\Ext^{i}_X(T,T)=0$ for all $i>0$, \item[(ii)] the endomorphism
algebra of T, $A=\Hom_X(T,T)$, has finite global homological
dimension, \item[(iii)] the direct summands of $T$ generate the
bounded derived category $D^b({\cO}_X)$.
\end{itemize}

\vskip 2mm The importance of tilting sheaves relies on the fact
that they can be characterized as those sheaves $T$ of
${\cO}_X$-modules such that the functors
\textbf{R}$\Hom_X(T,-):D^b(X)\longrightarrow D^b(A)$ and
$-\otimes_A^{\mbox{\textbf{L}}}T:D^b(A)\longrightarrow D^b(X)$
define mutually inverse equivalences of the bounded derived
categories of coherent sheaves on $X$ and of finitely generated
right $A$-modules, respectively. The existence of tilting sheaves also
plays also an important role in the problem of characterizing the
smooth projective varieties $X$ determined by its bounded category
 of coherent sheaves $D^b(X)$.
 For constructions of tilting
bundles and their relations to derived categories we refer to the
following papers: \cite{Ba}, \cite{Be}, \cite{Bo},
\cite{Ka} and \cite{Or}.

\vskip 2mm In this paper we will focus our attention on the
existence of tilting bundles on toric fibrations. The
search for tilting sheaves on a smooth projective variety $X$
splits naturally into two parts: First, we have to find the
so-called {\em strongly exceptional collection} of coherent
sheaves on $X$, $(F_0,F_1,\ldots ,F_n)$; and second we have to show that $F_0$, $F_1$,
$\ldots  ,$ $F_n$ generate the bounded derived category $D^b(X)$.

\begin{definition} A coherent sheaf $E$ on  a smooth projective variety  $X$ is called {\em exceptional} if
it is simple and $\Ext^{i}_{\cO _X}(E,E)=0$ for $i\ne 0$.  An ordered collection $(E_0,E_1,\ldots ,E_m)$  of coherent
sheaves on $X$ is an {\em exceptional collection} if each sheaf
$E_{i}$ is exceptional and $\Ext^i_{\cO _X}(E_{k},E_{j})=0$ for
$j<k$ and $i \geq 0$.
An exceptional collection $(E_0,E_1,\ldots ,E_m)$ is a {\em
strongly exceptional collection} if in addition $\Ext^{i}_{\cO
_X}(E_j,E_k)=0$ for $i\ge 1$ and  $j \leq k$. If an exceptional collection $(E_0,E_1,\cdots ,E_m)$ of coherent
sheaves on $X$ generates  $D^b(X)$, then it is called {\em full}.
\end{definition}
Thus each
full strongly exceptional collection defines a tilting sheaf
$T=\oplus_{i=0}^nF_{i}$ and, vice
versa,  each tilting bundle whose direct
summands are line bundles gives rise to a full strongly
exceptional collection.

\vspace{3mm}

The now classical result of Beilinson \cite{Be} states that $(\cO, \cO(1), \cdots , \cO(n))$ is a full strongly exceptional collection on $\PP^n$. The following problem can be considered by now a  natural and important question in Algebraic Geometry.

\begin{problem} Characterize smooth projective varieties which have a full strongly exceptional collection and investigate whether there is one consisting of line bundles.
\end{problem}

Note that not all smooth projective varieties have a full strongly exceptional collection of coherent sheaves. Indeed, the existence of a full strongly
exceptional collection $(E_0,E_1,\cdots,E_m)$ of coherent sheaves
on a smooth projective variety $X$ imposes a rather  strong
restriction on $X$, namely that the Grothendieck group
$K_0(X)=K_0({\cO}_X$-mod$)$ is isomorphic to $\ZZ^{m+1}$. Toric varieties fulfill a nice family of varieties that verify this condition on the Grothendieck group.
In \cite{CMZ}, \cite{CMPAMS}, \cite{CM} and \cite{CM1}, we constructed full strongly exceptional collections of line bundles on smooth toric varieties with a splitting fan, on smooth complete toric varieties with small Picard number and on $\PP^d$-bundles on certain toric varieties. The goal of this paper is to construct a full strongly exceptional collection consisting of line bundles
for toric fibrations and show the following result.

\begin{theorem} Let $X$ be a smooth complete toric variety. Assume that $X$ is a $F$-fibration over $Z$ and that $Z$ and $F$ are smooth toric varieties with a full strongly exceptional collection of line bundles. Then, $X$ has a full strongly exceptional collection made up of line bundles.
\end{theorem}

\vspace{3mm}
Two particular cases of this Theorem were proved by the authors in \cite{CMZ}. Indeed, in \cite{CMZ}, Theorem 4.17, we proved the case where $X$ is the trivial fibration over $Z$ with fiber $F$ (i.e. $X\cong Z\times F$), and in  \cite{CMPAMS}, Main Theorem we covered the case where $F\cong \PP^m$ for some integer $m$ .

\vspace{3mm}

Next we outline the structure of this paper. In section 2, we
recall the basic facts on toric varieties and toric fibrations.
Section 3 is the heart of this paper. First, we prove a key result
concerning the cohomology of line bundles on toric fibrations.
Then, in Theorem \ref{mainthm}, we prove that if $X$ is a
$F$-fibration over $Z$, and $Z$ and $F$  have both a full strongly
exceptional collection of line bundles then $X$ also has a full
strongly exceptional collection consisting of line bundles. We
divide the proof in two steps. First, we prove the existence of a
strongly exceptional collection on $X$.  Then, we prove that it is
full.

\vspace{3mm}

\noindent {\em Acknowledgment:} This work was initiated during a visit of the second author at the University of Barcelona, which she  thanks for its hospitality.


\section{An overview of Toric fibrations}

In this section we will introduce the notation and facts on toric
varieties and toric fibrations that we will use along this paper. We  refer to
\cite{Fu} and \cite{Oda} for more details.

\subsection{Toric varieties} Let $X$ be a smooth complete toric variety of dimension $n$ over
an algebraically closed field $K$ of characteristic zero. $X$ is defined by a
 fan $\Sigma:=\Sigma_X$ of strongly convex
polyhedral cones in $N\otimes _{\ZZ} \RR$ where $N$ is the lattice
$\ZZ^n.$ Let
$M:=\Hom_{\ZZ}(N,\ZZ)$ denote the dual lattice.

For any $0\le i \le n$, let $\Sigma (i):=\{
\sigma \in \Sigma \mid \text{ dim}(\sigma ) =i   \}$. In
particular, associated to any 1-dimensional cone $\sigma \in \Sigma(1)$ there
is a unique generator $v \in N$, called {\em ray generator}, such
that $\sigma \cap N=\ZZ_{\ge 0}\cdot v$. Denote by $G_X=\{v_{i} \mid i\in
J\}$ the set of ray generators of $X$.  There is a one-to-one correspondence between such ray
generators $\{v_{i} \mid i\in J\}$ and toric divisors $\{Z_{i}
\mid i\in J\}$ on $X$.

If $|G_X|=l$  then the Picard number of $X$ is
$\rho (X)=l-n$ and the anticanonical divisor $-K_X$ is given by
$-K_X=Z_1+\cdots+Z_l$.

\vskip 3mm Now we introduce  the notions of primitive collections
and primitive relations due to V.V. Batyrev \cite{Bat}.

\begin{definition}\label{primitrelat}
A set of toric divisors $\{Z_1,...,Z_k\}$ on $X$ is called a
{\em primitive set} if $Z_1\cap \cdot \cdot \cdot \cap Z_k=\emptyset $
but $Z_1\cap \cdot \cdot \cdot \cap \widehat{ Z_j}
\cap \cdot \cdot \cdot \cap Z_k\ne \emptyset$ for all $j$.
Equivalently, if $<v_1,...,v_k>\notin \Sigma$
but $<v_1,...,\widehat{v_j},...,v_k>\in \Sigma _X$ for all $j$.
$\cP=\{ v_1,...,v_k \}$ is called  a {\em primitive collection}.

If $\cP=\{ v_1,...,v_k \}
$ is a primitive collection, the element $v:=v_1+...+v_k$ lies
in the relative interior of a unique cone of $\Sigma_X $, say the cone
generated by $v_1',...,v_s'$ and
$v_1+...+v_k=a_1v_1'+...+a_sv_s'$ with $a_i>0$
is referred to as a  {\em primitive relation} $r(\cP)$ associated to $\cP$.
\end{definition}
 \subsection{Toric fibrations}
\begin{definition}
A {\it Toric Fiber-Bundle} is given by $(X,\phi, F, Z)$ where $X$ is a toric variety and $\phi$ a surjective morphism over a normal variety $Z$ such that every fiber is isomorphic to $F$.
\end{definition}

Observe that with these assumptions $Z$ and $F$ are necessarily toric and the map $\phi$ is equivariant. This implies that $\Sigma_F$ can be seen as a subfan of $\Sigma_X.$ The morphism $\phi$ is in fact  induced by an injective  map of fans $$i:\Sigma_F\hookrightarrow \Sigma_X.$$ Let $\dim(Z)=m, \dim(F)=d.$
The morphism $\phi$ being  locally trivial translates to the fact that for each $\sigma\in\Sigma_Z(m),$  $\phi^{-1}(U_{\sigma})\cong U_{\sigma}\times F,$ where $U_{\sigma}$ is the associated affine patch. This in turn implies that every $\sigma \in \Sigma_X(n)$ can be written as  $\sigma= \nu+\tau$, where $\nu\in\Sigma_F(d)$ and $\tau\cap\Sigma_F=\emptyset.$

This property characterizes in fact every fiber bundle:

\begin{proposition}[see for instance \cite{Ewa}; Theorem 6.7 or
\cite{Fu}; Section 2.4 for more details]\label{descripciofibration}
Let $X$ be an $n$-dimensional toric variety. $X$ has the structure of a toric fiber bundle  if and only if there is a linear subspace $H \subset N_{\RR}$
of dimension $d$ such that for every $n$-dimensional cone $\sigma \in \Sigma_X$, we have $\sigma= \nu+\tau$,
with $\nu \in \Sigma_X$, $\nu \subset H$, $\dim \nu = d$ and $\tau \cap H= \{0\}$.
\end{proposition}

In fact:
\begin{enumerate}
\item The set $\Sigma_F = \{\sigma \in \Sigma_X | \sigma \subset H  \}$ is the fan of a smooth, complete, $d$-dimensional toric variety $F$.
\item Denote by $H^ {\bot}$ the complementary space of $H$ in $N_{\RR}$ such that $N= (N \cap H ) \oplus (N \cap H^{\bot} )$
  and let $\pi: N_{\RR} \rightarrow H^{\bot}$ be the projection. Then the set
  \[ \Sigma_Z = \{\pi(\sigma) | \sigma \in \Sigma_X   \}\]
  is the fan of a complete, smooth $(n-d)$-dimensional toric variety $Z$.
\item  The projection $\pi$ induces an equivariant morphism $\phi: X \rightarrow Z$ such that, for every
  affine invariant open subset $U \subset Z$, ${\phi}^{-1}(U) \cong F \times U$ as toric varieties over $U$.
\end{enumerate}

We will use the name $F$-bundle for a fiber bundle with fiber $F.$

\begin{example}
(1) Let $F$ and $Z$ be smooth projective toric varieties. Then $X \cong F \times Z$ is the trivial fibration over $Z$
with fiber $F$.

(2) Let $X=\PP(\cE)$ be a smooth projective toric variety which is the projectivization
of a rank $r$ vector bundle $\cE$ on a smooth projective toric variety $Z$. Then, $X$ is a  $\PP^{r-1}$-fibration over $Z$.

\end{example}

Let $X$ be an $F$-fibration over $Z$ as in Proposition \ref{descripciofibration}. Then $G_F=G_X \cap H$ and the primitive
collections of $\Sigma_X$ contained in $H$ are exactly the primitive collections of $\Sigma_F$. These primitive
collections have the same primitive relations in $\Sigma_X$ and in $\Sigma_F$. On the other hand, the projection
$\pi: N_{\RR} \rightarrow H^{\bot}$ induces a bijection between $G_X \backslash G_F$ and $G_Z$. Under this identification, the
primitive collections of $\Sigma_X$ not contained in $H$ are exactly the primitive collections of $\Sigma_Z$.
If we denote by $\overline{x}$ the image of $x \in N$ by $\pi$, then a primitive collection in $\Sigma_Z$
of the form \[ \overline{x}_1+ \cdots + \overline{x}_h -(a_1 \overline{y}_1 + \cdots + a_k \overline{y}_k)=0 \]
lifts in $\Sigma_X$   to a primitive relation of the form
 \[ x_1+ \cdots + x_h -(a_1 y_1 + \cdots + a_k y_k+b_1 z_1+ \cdots +b_lz_l)=0, \]
 with $l \geq 0$ and $z_i \in \Sigma_F$ for $1 \leq i \leq l$. In particular, the fibration is trivial, namely
 $X \cong F \times Z$, if and only if all primitive relations in $\Sigma_Z$ remain unchanged when lifted
 to $\Sigma_X$.

\vspace{3mm}

For any smooth projective toric variety $X$, we denote by $P_X(t)$
its Poincar\'e polynomial. It is well known that the topological
Euler characteristic of $X$, $\chi(X)$ verifies
\[ \chi(X)=P_X(-1)\]
and $\chi(X)$ coincides with the number of maximal cones of $X$,
that is, with the rank of the Grothendieck group $K_0(X)$ of $X$.
On the other hand, if  $X$ is an $F$-fibration over $Z$ we have
(\cite{Fu};Pg. 92-93):

\[ P_{X}(t)=P_{F}(t)\cdot P_{Z}(t).\]

Thus putting altogether we deduce that:

\begin{equation} \label{rankK0} \rank(K_0(X))=  P_{X}(-1)=P_{F}(-1) \cdot P_{Z}(-1)
=\rank(K_0(F))\cdot \rank(K_0(Z)). \end{equation}

\vspace{3mm}

Let us now describe the toric $F$-fibration $\phi: X \rightarrow Z$ in terms of the smooth toric varieties $F$ and $Z$.

\vspace{3mm}

Let  $\dim(Z)=m,\dim(F)=d$ and as in the previous section we
denote by  $\Sigma_Z \subset  \RR^m,$ and $\Sigma_F\subset \RR^d$
the defining fans.

Consider $\sigma= \langle  v_1, \ldots, v_d
\rangle,$ a maximal cone of $F$. Since $F$ is a smooth toric
variety,  we can  take
${v}_1, \ldots, {v}_d$
to be a $\ZZ$-basis of $\ZZ^d.$
Let $G_F=\{v_1, \ldots,
{v}_d,{v}_{d+1}, \ldots, {v}_{r} \}$
be the set of ray generators of $F$ expressed in this $\ZZ$-basis and
denote by $\tilde{F}_i$ the corresponding toric divisors.  We fix
\[  \langle \tilde{F}_{d+1}, \tilde{F}_{d+2}, \cdots, \tilde{F}_{r} \rangle \]
to be a
basis of $\Pic(F).$ In this basis we have the following linear relations
\begin{equation} \label{pic1} \tilde{F}_i = \sum_{j=d+1}^{r} \delta_i^j \tilde{F}_j, \quad \quad 1 \leq i \leq d, \quad \mbox{ for some } \delta_i^j \in \ZZ. \end{equation}
Consider $\tau= \langle  e_1, \ldots, e_m \rangle,$ a maximal cone
of $Z$. Again, since $Z$ is a smooth toric variety,  we can  take
$e_{1}, \ldots, e_{m}$ to be a $\ZZ$-basis of $\ZZ^m$. Let
$G_Z=\{e_1, \ldots, e_m,e_{m+1}, \ldots, e_{t} \}$ be the set of
ray generators of $Z$ expressed in this $\ZZ$-basis and denote by
$\tilde{Z}_i$ the corresponding toric divisors.  We fix
\[  \langle \tilde{Z}_{m+1}, \tilde{Z}_{m+2}, \cdots, \tilde{Z}_{t} \rangle \]
to be a
basis of $\Pic(Z).$ In this basis we have the following linear relations
\begin{equation} \label{pic2}  \tilde{Z}_i = \sum_{j=m+1}^{t} \alpha_i^j \tilde{Z}_j, \quad \quad 1 \leq i \leq m, \quad \mbox{ for some } \alpha_i^j \in \ZZ. \end{equation}
By the way the fan of a toric-bundle is constructed one sees that in the $\ZZ$-basis of $\ZZ^{d+m}$
$$ \{({v}_1,0), \cdots, ({v}_d,0),(0,e_{1}), \ldots, (0,e_{m})\}$$ the fan of $X$ has ray generators
$$G_X=\{\overline{v_1}, \ldots,
\overline{v_d},\overline{v_{d+1}}, \ldots, \overline{v_{r}},\overline{e_1}, \ldots,
\overline{ e_m},\overline{e_{m+1}}, \ldots, \overline{e_{t}}\}$$
with for some $\gamma_k^j \in \ZZ$,
\begin{equation} \label{invariants} \begin{array}{l}
\overline{v_i}:=(v_i,0), \quad 1 \leq i \leq r \\
\overline{e_i}:=(0,e_i), \quad 1 \leq i \leq m \\
\overline{e_j}:= (\sum_{k=1}^{d} \gamma_k^j v_k, e_j), \quad m+1 \leq j \leq t;
\end{array}\end{equation}
and the projection $\pi:\ZZ^{d+m}\to\ZZ^m$ is so that $\pi(\overline{v_i})=0$ for $i=1,\ldots r$ and $\pi(\overline{e_j})=e_j$ for $j=1,\ldots t.$
It follows that $${\rm Pic}(X)=<F_{d+1},\ldots,F_r,Z_{m+1},\ldots,Z_t>$$
where the $F_i$ are the toric divisors associated to the $\overline{v_i}$ and the $Z_j$ are the divisors associated to the $\overline{e_j}.$

According to (\ref{pic1}), (\ref{pic2}) and (\ref{invariants}), we have the following linear relations
\begin{equation}
\label{picX}
\begin{array}{l}
Z_i= \sum_{k=m+1}^{t} \alpha_i^k Z_k, \quad 1 \leq i \leq m \\
F_i= \sum_{j=d+1}^{r} \delta_i^j F_j - \sum_{j=m+1}^{t} \gamma_i^j Z_j, \quad 1 \leq i \leq d.
\end{array}
\end{equation}
\vspace{3mm}

We end the section with the following remark that will be used frequently during this work.

\begin{remark}
\label{isoclau}
From the construction of the $F$-fibration $\phi: X \rightarrow Z$, it follows that
\[ \cO_X(F_i)|_F=\cO_F(\tilde{F_i}) \text{ and } \cO_X(Z_i)=\phi^* \cO_Z(\tilde{Z_i}) \]
and given line bundles $L=\cO_F(\sum_{i=d+1}^{r}a_i \tilde{F_i})$ on $F$ and $H=\cO_Z(\sum_{j=m+1}^{t}b_j \tilde{Z_j})$ on $Z$, the line bundle
\[ \cL= \cO_X(\sum_{i=d+1}^{r}a_i {F_i} + \sum_{j=m+1}^{t}b_j {Z_j}) \cong
\cO_X(\sum_{i=d+1}^{r}a_i {F_i}) \otimes \phi^* \cO_Z(\sum_{j=m+1}^{t}b_j \tilde{Z_j})\]
on $X$ is so that
\[ \cL|_F \cong  \cO_F(\sum_{i=d+1}^{r}a_i \tilde{F_i})=L. \]
\end{remark}

\section{Derived category of toric fibrations}

This section contains the main new results of this work. Its
  goal is to give a structure theorem for the derived category $D^b(X)$ of a toric $F$-fibration $X$ over $Z$. This will be achieved  constructing a full strongly exceptional collection
of line bundles on $X$ being $X$ a toric fibration over $Z$ with fibers a toric variety $F$ such that the  derived categories $D^b(Z)$ and $D^b(F)$ have an orthogonal basis consisting  of line bundles.

We start by recalling the
notions of exceptional sheaves, exceptional collections of
sheaves, strongly exceptional collections of sheaves and full
strongly exceptional collections of sheaves as well as the cohomology of line bundles on toric
$F$-fibrations needed in the sequel.

\begin{definition}\label{exceptcoll}
Let $X$ be a smooth projective variety.

(i) A  coherent sheaf $E$ on  $X$ is {\em exceptional} if $\Hom
(E,E)=K $ and $\Ext^{i}_{\cO _X}(E,E)=0$ for $i>0$,

(ii) An ordered collection $(E_0,E_1,\ldots ,E_m)$  of coherent
sheaves on $X$ is an {\em exceptional collection} if each sheaf
$E_{i}$ is exceptional and $\Ext^i_{\cO _X}(E_{k},E_{j})=0$ for
$j<k$ and $i \geq 0$.

(iii) An exceptional collection $(E_0,E_1,\ldots ,E_m)$ is a {\em
strongly exceptional collection} if in addition $\Ext^{i}_{\cO
_X}(E_j,E_k)=0$ for $i\ge 1$ and  $j \leq k$.

(iv) An ordered collection $(E_0,E_1,\cdots ,E_m)$ of coherent
sheaves on $X$ is a {\em full (strongly) exceptional collection}
if it is a (strongly) exceptional collection  and $E_0$, $E_1$,
$\cdots $ , $E_m$ generate the bounded derived category $D^b(X)$.
\end{definition}

\begin{remark} \label{length} As mentioned in the Introduction, the existence of a full strongly
exceptional collection $(E_0,E_1,\cdots,E_m)$ of coherent sheaves
on a smooth projective variety $X$ imposes a rather  strong
restriction on $X$, namely that the Grothendieck group
$K_0(X)=K_0({\cO}_X$-mod$)$ is isomorphic to $\ZZ^{m+1}$.
\end{remark}

Let us illustrate the above definitions with  examples:

\begin{example} \label{prihirse}

(1) Suppose that the vectors $v_0, \cdots, v_n$ generate a lattice $N$ of rank $n$ and
$v_0+\cdots+v_n=0$. Let $\Sigma$  be the fan whose rational cones are generated by any proper subsets of the
vectors $v_0, \cdots, v_n$. It is well known that  the toric variety associated to the fan $\Sigma$ is  $\PP^n$.
Denote by $Z$ the toric divisor associated to $v_0$. Then the collection of line bundles
$(\cO_{\PP^n}, \cO_{\PP^n}(Z), \cdots, \cO_{\PP^n}(nZ))$ is a full strongly exceptional collection of line bundles on
$\PP^n$.

\vskip 2mm (2) Let $X_1$ and $X_2$ be two smooth projective
varieties and let $(F_0^{i},F_1^{i},\ldots ,F_{n_{i}}^{i})$ be   a
full strongly exceptional collection of locally free sheaves on
$X_i$, $i=1,2$.  Then, $$(F_0^{1}\boxtimes
F_0^{2},F_1^{1}\boxtimes F_0^{2},\ldots ,F_{n_1}^{1}\boxtimes
F_0^{2},  \ldots , F_0^{1}\boxtimes F_{n_2}^{2},F_1^{1}\boxtimes
F_{n_2}^{2},\ldots ,F_{n_1}^{1}\boxtimes F_{n_2}^{2})$$ is a full
strongly exceptional collection of locally free sheaves on $X_1
\times X_2$ (see \cite{CMZ}; Proposition 4.16). In particular, if $Z_n$ is a toric divisor that generates $\Pic(\PP^n)$ and
$Z_m$ is a toric divisor that generates $\Pic(\PP^m)$, then
the collection of line bundles
$$( \cO_{\PP^n} \boxtimes \cO_{\PP^m},
\cO_{\PP^n}(Z_n) \boxtimes \cO_{\PP^m},
\cdots, \cO_{\PP^n}(nZ_n) \boxtimes \cO_{\PP^m},
\cdots,  \cO_{\PP^n} \boxtimes \cO_{\PP^m}(mZ_m), $$ $$ \hspace{20mm} \cO_{\PP^n}(Z_n)
\boxtimes \cO_{\PP^m}(mZ_m),\cdots, \cO_{\PP^n}(nZ_n) \boxtimes
\cO_{\PP^m}(mZ_m)) $$ is a full strongly exceptional collection of
line bundles on $\PP^n \times \PP^m$.

\vskip 2mm (3) Let $X$ be a smooth complete toric variety which is
the projectivization of a rank $r$ vector bundle $\cE$ over a
smooth complete toric variety $Z$ which has a full strongly
exceptional collection of locally free sheaves. Then, $X$ also has
a full strongly exceptional collection of locally free sheaves
(See \cite{CMZ}; Proposition 4.9).

\end{example}

\vspace{4mm}

As we said in the introduction $D^b(X)$ is an important algebraic invariant of a smooth projective variety but very little is known about the structure of $D^b(X)$. In particular, whether  $D^b(X)$ is freely and finitely
generated and, hence, we are lead to consider the following
problem

\begin{problem}\label{prob1}
Characterize smooth projective varieties $X$ which have a full
strongly exceptional collection of coherent sheaves and, even
more, if there is one made up of line bundles.
\end{problem}

This problem is far from being solved and in this paper we will
restrict our attention to the particular case of toric varieties since they
admit a combinatorial description which allows
many invariants to be expressed
in terms of combinatorial data.

In \cite{Kaw}, Kawamata proved that the derived category of a
smooth complete toric variety has a full exceptional collection of
objects. In his collection the objects are sheaves rather than
line bundles and the collection is only exceptional and not
strongly exceptional. In the toric context,  there are a lot of
contributions to the above problem. For instance, it turns out
that a full strongly exceptional collection made up of line
bundles exists on  projective spaces (\cite{Be}), multiprojective
spaces (\cite{CMZ}; Proposition 4.16), smooth complete toric
varieties with Picard number $\le 2$ (\cite{CMZ}; Corollary 4.13), smooth complete toric varieties with a splitting fan
(\cite{CMZ}; theorem 4.12) and some smooth complete toric varieties with Picard Number 3 (\cite{MM}). Nevertheless some restrictions are
required because in \cite{HP}, Hille and Perling
constructed an example of smooth non Fano toric surface which does
not have a full strongly exceptional collection made up of line
bundles. Recently in \cite{Mateusz} M. Michalek has also proved that there is no a strongly exceptional collection of line bundles on the blow up of $\PP^n$, $n$ big enough,  at two $T$-invariant points.

\vspace{3mm}

Our goal is to describe the derived category of a toric fibration. This will be achieved after an accurate description of acyclic line bundles on this varieties.

  \vspace{3mm}

\begin{definition} Let $X$ be a  smooth complete toric variety. A
line bundle $\cL$ on $X$ is said to be {\em acyclic} if
$H^{i}(X,\cL)=0$ for every $i\ge 1$.
\end{definition}

\vspace{3mm}

Now we are going to provide a description of cohomology of line bundles $\cL$ on an $n$-dimensional
smooth complete toric variety  $X$. To this end we need to fix some
notation.

\begin{notation} Let $X$ be a smooth complete toric variety and denote by $T_1, \cdots, T_s$ the toric divisors on $X$.
For every $ \textbf{r}=(r_i)_{i=1}^{s} \in \ZZ^s$ we denote by $\Supp(\textbf{r})$ the simplicial complex on
$s$ vertices, which encodes all the cones $\sigma$ of $\Sigma_X$ for which all $i$ with $v_i$ in $\sigma$
satisfy $r_i \geq 0$. As usual $H_q(\Supp(\textbf{r}), K)$ denotes the $q$-th reduced
homology group
of the simplicial complex $\Supp(\textbf{r})$.
\end{notation}

\vspace{3mm}

The following lemma which is a
direct consequence of well known facts concerning homology of
simplicial complexes will be very useful in the sequel.

\vspace{3mm}

\begin{lemma}\label{key1}
\label{kunneth} Assume that $
\textbf{r}=(r_i)_{i=1}^{r+t} \in \ZZ^{r+t}$ decomposes into $\textbf{r} =(\textbf{r}_1, \textbf{r}_2) \in \ZZ^r \times \ZZ^t$. Then the following K\"{u}nneth formula
holds:
\[ H_i(\Supp(\textbf{r}), K) \cong \bigoplus_{p+q=i} H_p(\Supp(\textbf{r}_1), K) {\cdot} H_q(\Supp(\textbf{r}_2), K). \]
\end{lemma}

 Given a line bundle $\cL$ on $X$,
we will say that $\textbf{r}=(r_i)_{i=1}^{s} \in \ZZ^s$ represents $\cL$ whenever  $\cL\cong  \cO_{X}(r_1 T_1+r_2 T_2+\cdots+r_sT_s)$. We have

\vspace{3mm}

\begin{proposition}\label{cohomology}
With the above notation, we have
$$H^p(X, \cL)\cong \bigoplus _ {\textbf{r}}H_{\rank{N}-p}(\Supp(\textbf{r}), K)$$
where  the sum is taken over all the representations $\textbf{r}=(r_i)_{i=1}^{s} \in \ZZ^s$ of $\cL$.
\end{proposition}
\begin{proof} See \cite{bh}; Proposition 4.1.
\end{proof}

\vspace{3mm} As a consequence we obtain

\begin{corollary} \label{h0}
With the above notation, $H^0(X, \cL)$ is determined only by  $
\textbf{r}=(r_i)_{i=1}^{s} \in \ZZ^s$ such that
$\Supp(\textbf{r})$ is the entire fan $\Sigma_X.$ Equivalently  $H^0(X, \cL)$ is determined only by $\textbf{r}=(r_i)_{i=1}^{s} \in \ZZ^s$  such that
$\cL\cong  \cO_{X}(r_1 T_1+r_2 T_2+\cdots+r_sT_s)$ with all
$r_{i}\ge 0$. On the other hand, $H^{\dim X}(X, \cL)$ is determined only
by
 $ \textbf{r}=(r_i)_{i=1}^{s} \in \ZZ^s$
such that $\Supp(\textbf{r})$ is empty, i.e. by  $ \textbf{r}=(r_i)_{i=1}^{s} \in \ZZ^s$ such that $\cL\cong
\cO_{X}(r_1 T_1+r_2 T_2+\cdots+r_sT_s)$ with all $r_{i}\le
-1$.
\end{corollary}

Our next goal is to somehow control all the representations of a line bundle on a toric fiber bundle. We have

\begin{lemma} \label{key2} Let $\phi:X \longrightarrow Z$ be a toric $F$-fibration and keep the notation introduced in \S 2. Assume that $\textbf{r}=(a_1, \cdots, a_r,b_1, \cdots , b_t)$ is a representation of
$$\cD =\cO_X(\sum _{i=d+1}^r\alpha _{i}F_{i}+\sum _{j=m+1}^t\beta _jZ_j)=
\cO_X(\sum _{i=d+1}^r\alpha _{i}F_{i}) \otimes \phi^* \cO_Z(\sum _{j=m+1}^t\beta _j\tilde{Z}_j). $$ Then,  $
\textbf{r}=(\textbf{r}_1,\textbf{r}_2)$ where  $\textbf{r}_1= (a_1, \cdots, a_r)$ is a representation of
 $ \cO_F(\sum _{i=d+1}^r\alpha _{i}\tilde{F}_{i})$ and
 $\textbf{r}_2=  (b_1, \cdots , b_t)$  is a representation of
$\cO_Z(\sum _{j=m+1}^t\beta _j\tilde{Z}_j+D^{\gamma}_{a_1,\cdots , a_r})$
 with
$$D^{\gamma}_{a_1, \cdots, a_r}:=(\sum _{i=1}^{d}a_i\gamma _{i}^{m+1})\tilde{Z}_{m+1}+\cdots +(\sum _{i=1}^{d}a_i\gamma _{i}^{t})\tilde{Z}_{t}.$$
\end{lemma}
\begin{proof} Using the relations (\ref{picX}) together with the fact that  $\textbf{r}=(a_1, \cdots, a_r,b_1, \cdots , b_t)$ is a representation of $\cD$, we obtain
$$\begin{array}{ccl} \sum _{i=d+1}^r\alpha _{i}F_{i}+\sum _{j=m+1}^t\beta _jZ_j & = & \sum _{i=1}^r a_{i}F_{i} + \sum_{j=1}^{t}b_jZ_j \\  \\ & = &
\sum _{j=1}^{m} b_{j}(\sum_{k=m+1}^{t} \alpha_j^kZ_k ) + \sum _{j=m+1}^t b_{j}Z_{j} \\ & & +\sum _{i=1}^{d}  a_{i}(\sum_{j=d+1}^{r} \delta_i^jF_j - \sum_{j=m+1}^{t} \gamma_i^jZ_j)\\ & & + \sum _{i=d+1}^r a_{i}F_{i}.\end{array} $$

Therefore, we have
$$ \begin{array}{ccl}  \alpha_{d+1} & = & \sum _{i=1}^{d} a_{i}\delta _{i}^{d+1}+ a_{d+1} \\ \vdots  & & \\
\alpha_{r} & = & \sum _{i=1}^{d} a_{i}\delta _{i}^{r}+ a_{r} \\
\beta_{m+1} & = & \sum _{j=1}^m b_{j}\alpha_{j}^{m+1}+ b_{m+1} - \sum _{i=1}^{d} a_{i} \gamma _{i}^{m+1} \\ \vdots  & & \\  \beta_{t} & = & \sum _{j=1}^m b_{j}\alpha _{j}^{t}+ b_{t} - \sum _{i=1}^{d} a_{i} \gamma _{i}^{t} \\
\end{array}$$ which gives us that  $\textbf{r}_1=  (a_1, \cdots, a_r)$  is a representation of
  $ \cO_F(\sum _{i=d+1}^r\alpha _{i}\tilde{F}_{i})$ and
 $\textbf{r}_2=  (b_1, \cdots , b_t)$  is a representation of
$\cO_Z(\sum _{j=m+1}^t\beta _j\tilde{Z}_j+D^{\gamma}_{a_1,\cdots , a_r})$
 with
$$D^{\gamma}_{a_1, \cdots, a_r}=(\sum _{i=1}^{d}a_i\gamma _{i}^{m+1})\tilde{Z}_{m+1}+\cdots +(\sum _{i=1}^{d}a_i\gamma _{i}^{t})\tilde{Z}_{t}.$$
\end{proof}

\begin{remark} It is important to point out that the divisor  $D^{\gamma}_{a_1,\cdots , a_r}$ does not depend on
$\phi^* \cO_Z(\sum _{j=m+1}^t\beta _j\tilde{Z}_j)$. In fact, the $\gamma _{i}^{j}$ are invariants of the fibration and $(a_1,\cdots , a_r)$ is a representation of $ \cO_F(\sum _{i=d+1}^r\alpha _{i}\tilde{F}_{i})$.
\end{remark}

\vspace{3mm}

\begin{corollary}\label{acyclic}With the above notation,  let $\sum _{i=d+1}^r\beta _i\tilde{F}_i$ be a non-effective divisor on $F$ such that  the corresponding line bundle $\cO_F(\sum _{i=d+1}^r\beta _i\tilde{F}_i)$ is acyclic. Then, for any line bundle $\cO_Z(\sum _{i=m+1}^t\alpha _{i}\tilde{Z}_{i})$ on $Z$,
$$\cO_X(\sum _{i=d+1}^r\beta _iF_i ) \otimes \phi^* \cO_Z(\sum _{i=m+1}^t\alpha _{i}\tilde{Z}_{i})$$ is an acyclic line bundle on $X$.

\end{corollary}
\begin{proof} It follows from Lemma \ref{key1}, Proposition \ref{cohomology} and Lemma \ref{key2}.
\end{proof}

\begin{notation} From now on, given a line bundle $\tilde{\cL}=\cO_F(\sum_{i=d+1}^{r}a_i\tilde{F}_i)$ on $F$ we denote by ${\cL}=\cO_X(\sum_{i=d+1}^{r}a_i{F}_i)$ the line bundle on $X$ so that $\cL|_F=\tilde{\cL}$.
\end{notation}

We are now ready to state the main result of this work.

\vspace{3mm}

\begin{theorem} \label{mainthm} Let $\phi: X \rightarrow Z$ be a toric $F$-fibration over $Z$. Assume that $(\tilde{\cL_1}, \cdots , \tilde{\cL_u})$ is a full  strongly exceptional
collection of line bundles on $D^b(F)$ and that  $(\cE_1, \cdots , \cE_v)$ is a  full strongly exceptional collection of line bundles on  $D^b(Z)$. Then, there is a line bundle $\cO_Z(D)$ on $Z$ such that
the ordered sequence
$$\{\{ \phi^*(\cE_j \otimes \cO_Z(D))\otimes \cL_1\}_{ 1\leq j\leq v}, \{\phi^*(\cE_j \otimes \cO_Z(2D))\otimes \cL_2\}_{ 1\leq j\leq v}, \ldots, \{ \phi^*(\cE_j \otimes \cO_Z(uD))\otimes \cL_u\}_{ 1\leq j\leq v}\}$$
is a full strongly exceptional sequence of line bundles on $X.$
\end{theorem}
\begin{proof} We have to prove the existence a line bundle $\cO_Z(D)$ on $Z$ such that the above collection is full and satisfies the following cohomological conditions:

(a) $\Ext^k( \phi^*(\cE_i \otimes \cO_Z(jD))\otimes \cL_j, \phi^*(\cE_p \otimes \cO_Z(qD)) \otimes \cL_q) =0$ for $k > 0$ and $q \geq j$, and

(b) $\Ext^k(\phi^*(\cE_i \otimes \cO_Z(jD)) \otimes \cL_j , \phi^*(\cE_p \otimes \cO_Z(qD)) \otimes \cL_q)=0$ for $k \geq 0$ and $q < j$.

\noindent Let us check $(a)$. We distinguish two cases. First of all assume that $q >j$. In that case,
$$\Ext^k( \phi^*(\cE_i \otimes \cO_Z(jD))\otimes \cL_j, \phi^*(\cE_p \otimes \cO_Z(qD)) \otimes \cL_q)
=H^k(\phi^*(\cE_p \otimes \cE_i^{\vee}  \otimes \cO_Z((q-j)D)) \otimes \cL_q \otimes \cL_j^{\vee}).$$
If this cohomology group is non-zero, by Proposition \ref{cohomology},  there exists $\textbf{r}=(a_1, \cdots,a_r,b_1, \cdots, b_t)$ a representation of $\phi^*(\cE_p \otimes \cE_i^{\vee}  \otimes \cO_Z((q-j)D)) \otimes \cL_q \otimes \cL_j^{\vee}$ such that
$$ H_{n-k}(\Supp(\textbf{r}), K)  \neq 0, \quad n-k < n. $$

By Lemma \ref{key2} and Lemma \ref{key1}, this means that there exists $\textbf{r}_1=(a_1, \cdots, a_r)$ a representation of $ \tilde{\cL_q} \otimes \tilde{\cL_j}^{\vee}$ with $q >j$ and $\textbf{r}_2=(b_1, \cdots, b_t)$ a representation of $$ \cE_p \otimes \cE_i^{\vee}  \otimes \cO_Z((q-j)D) \otimes \cO_Z(D^{\gamma}_{a_1, \cdots,a_r})$$ 
such that for some $f$ and $g$ with $f+g=n-k$,
$$ H_{f}(\Supp(\textbf{r}_1), K)  \neq 0 \quad \mbox{and} \quad
 H_{g}(\Supp(\textbf{r}_2), K)  \neq 0. $$
 Since $(\tilde{\cL_1}, \cdots , \tilde{\cL_u})$ is a  strongly exceptional collection of line bundles in  $D^b(F)$,
 \[ H^t (\tilde{\cL_q }\otimes \tilde{\cL_j}^{\vee})=0 \quad \mbox{for any} \quad t>0.\]
 Hence,  the representation $\textbf{r}_1$ must contribute to $H^0 (\tilde{\cL_q} \otimes \tilde{\cL_j}^{\vee})$ and by Corollary \ref{h0} this implies that
 $a_j \geq 0$ for $ 1 \leq j \leq r$. Moreover, since $H^0 (\tilde{\cL_q} \otimes \tilde{\cL_j}^{\vee})$ has finite dimension, there is only a finite number of such representations $\textbf{r}_1=(a_1, \cdots, a_r)$ of $ \tilde{\cL_q} \otimes \tilde{\cL_j}^{\vee}$. Thus, since the coefficients $\gamma_i^j$ are invariants of the variety and $(\cE_1, \dots, \cE_v)$ is a finite collection, the divisors $D^{\gamma}_{a_1, \cdots,a_r}$ move in a finite set. Therefore, since $q-j>0$ we can take a big enough line bundle $\cO_Z(D)$ on $Z$ such that 
$\cE_p \otimes \cE_i^{\vee}  \otimes \cO_Z((q-j)D) \otimes \cO_Z(D^{\gamma}_{a_1, \cdots,a_r})$ is globally generated. Hence 
   $H^l( \cE_p \otimes \cE_i^{\vee}  \otimes \cO_Z((q-j)D) \otimes \cO_Z(D^{\gamma}_{a_1, \cdots,a_r}))=0$ for any $l>0$ which will contradict the fact that
 \[  H_{g}(\Supp(\textbf{r}_2), K)  \neq 0 \]
 being $\textbf{r}_2$  a representation of $ \cE_p \otimes \cE_i^{\vee}  \otimes \cO_Z((q-j)D) \otimes \cO_Z(D^{\gamma}_{a_1, \cdots,a_r})$.

 In case, $q=j$, we have
 $$\Ext^k( \phi^*(\cE_i \otimes \cO_Z(jD))\otimes \cL_j, \phi^*(\cE_p \otimes \cO_Z(jD)) \otimes \cL_j)
=H^k(\phi^*(\cE_p \otimes \cE_i^{\vee} ))=0 \quad \mbox{for all} \quad k>0$$
 since by assumption $(\cE_1, \cdots , \cE_v)$ is a full  strongly exceptional
collection of line bundles in $D^b(Z)$ and in this case, by Lemma \ref{key1}, the cohomology on $X$ coincides  with the cohomology on $Z$.

In case (b), if there exist $k \geq 0$ and $q < j$ such that
$$\Ext^k(\phi^*(\cE_i \otimes \cO_Z(jD)) \otimes \cL_j , \phi^*(\cE_p \otimes \cO_Z(qD)) \otimes \cL_q)=$$
$$ H^k(\phi^*(\cE_p \otimes \cE_i^{\vee}  \otimes \cO_Z((q-j)D)) \otimes \cL_q \otimes \cL_j^{\vee})\neq 0, $$
then by Proposition \ref{cohomology}, there exists $\textbf{r}=(a_1, \cdots, a_r,b_1, \cdots, b_t )$ a representation of the line bundle $\phi^*(\cE_p \otimes \cE_i^{\vee}  \otimes \cO_Z((q-j)D)) \otimes \cL_q \otimes \cL_j^{\vee}$ such that
$$ H_{n-k}(\Supp(\textbf{r}), K)  \neq 0. $$
By Lemma \ref{key2} and Lemma \ref{key1}, this means that there exists $\textbf{r}_1=(a_1, \cdots, a_r)$ a representation of $ \tilde{\cL_q} \otimes \tilde{\cL_j}^{\vee}$ with $q < j$  such that for some $f$,
$$ H_{f}(\Supp(\textbf{r}_1), K)  \neq 0. $$
But this contradicts the fact that  $(\tilde{\cL_1}, \cdots , \tilde{\cL_u})$ is a full strongly exceptional collection of line bundles in  $D^b(F)$.

Putting altogether we get that indeed there exists an ample line bundle $\cO_Z(D)$ on $Z$ such that the collection is strongly exceptional.

\vspace{3mm}
Let us prove that it is full. To this end, observe that since by assumption $(\cE_1, \cdots,\cE_v )$ is a full strongly exceptional collection on $Z$, for any $j>0$
\[ D^b(Z)= \langle \cE_1 \otimes \cO_Z(jD), \cdots,\cE_v \otimes \cO_Z(jD) \rangle.\]
Indeed, since clearly  $( \cE_1 \otimes \cO_Z(jD), \cdots,\cE_v \otimes \cO_Z(jD))$ is also an exceptional collection, to see that it is full it is enough (\cite{BO}) to see that
\begin{equation} \label{orto} 0= {}^{\bot} \langle \cE_1 \otimes \cO_Z(jD), \cdots,\cE_v \otimes \cO_Z(jD) \rangle
:= \{ F \in D^b(Z)| Ext^{\bullet}(F,\cE_i \otimes \cO_Z(jD) )=0 \}. \end{equation}
By assumption, $(\cE_1, \cdots,\cE_v )$ is a full strongly exceptional collection. Hence
$$ {}^{\bot} \langle \cE_1, \cdots,\cE_v  \rangle := \{ G \in D^b(Z)| Ext^{\bullet}(G,\cE_i )=0\} =0 $$
which is equivalent to say that the equality in (\ref{orto}) holds.

Therefore, it suffices  to prove that
\[D^b(X)=\langle \phi^*D^b(Z) \otimes \cL_1, \cdots, \phi^*D^b(Z) \otimes \cL_u \rangle \]
where $\phi^*D^b(Z) \otimes \cL_k$ denotes the full triangulated subcategory in $D^b(X)$ generated by the objects
$\{ \phi^*(\cE_j)\otimes \cL_k \}_{ 1\leq j\leq v}$.

To this end, it is enough to check that  $ \langle \phi^*D^b(Z) \otimes \cL_1, \cdots, \phi^*D^b(Z) \otimes \cL_u \rangle$ contains all the objects $\cO_x$, $x \in X$ since the set $\{\cO_x| x \in X \}$ is a spanning class for $D^b(X)$. Since every point $x \in X$ belongs to some fiber $F_z:= \phi^{-1}(z)$, $z \in Z$, and by assumption
$(\tilde{\cL_1}, \cdots , \tilde{\cL_u})$ generates $D^b(F)$,
\[ \cO_x \in \langle \tilde{\cL_1}, \cdots , \tilde{\cL_u} \rangle  .\]
So, since $\phi^*(\cO_z)=\cO_{F_z}$, the sheaf $\cO_x$ belongs to
$$\langle \phi^*D^b(Z) \otimes \cL_1, \cdots, \phi^*D^b(Z) \otimes \cL_u \rangle.$$
Putting altogether we get that
$$\{\{ \phi^*(\cE_j \otimes \cO_Z(D))\otimes \cL_1\}_{ 1\leq j\leq v}, \{\phi^*(\cE_j \otimes \cO_Z(2D))\otimes \cL_2\}_{ 1\leq j\leq v}, \ldots, \{ \phi^*(\cE_j \otimes \cO_Z(uD))\otimes \cL_u\}_{ 1\leq j\leq v}\}$$
generates $D^b(X)$.
  \end{proof}

\vspace{3mm}

\begin{remark} 
Let $\phi: X \rightarrow \PP^n$ be a toric $F$-fibration over $\PP^n$. Assume that $(\tilde{\cL_1}, \cdots , \tilde{\cL_u})$ is a full  strongly exceptional
collection of line bundles on $D^b(F)$ and consider  $(\cO_{\PP^n}, \cdots , \cO_{\PP^n}(n))$ a  full strongly exceptional collection of line bundles on  $D^b(\PP^n)$. In that particular case, following the same ideas as in the proof of Theorem \ref{mainthm}, we get that the ordered sequence
$$\{\{ \phi^*(\cO_{\PP^n})\otimes \cL_1\}_{ 1\leq j\leq v}, \{\phi^*(\cO_{\PP^n}(1))\otimes \cL_2\}_{ 1\leq j\leq v}, \ldots, \{ \phi^*(\cO_{\PP^n}(n) )\otimes \cL_u\}_{ 1\leq j\leq v}\}$$
is a full strongly exceptional sequence of line bundles on $X$.
\end{remark}



\begin{thebibliography}{999}

\bibitem{Ba} D. Baer, {\em  Tilting sheaves}, Manusc. Math. {\bf 60} (1988), 323a-347.



\bibitem{Bat} V.V.  Batyrev,
{\em On the classification of smooth projective
toric varieties}, Tohoku Math. J. {\bf 43} (1991), 569-585.




\bibitem{Be} A.A. Beilinson, {\em Coherent sheaves on $\PP^n$ and Problems of Linear Algebra},
Funkt. Anal. Appl. {\bf 12} (1979), 214-216.


 \bibitem{Bo} A.I. Bondal, {\em Representation of associative algebras
  and coherent sheaves}, Math. USSR Izvestiya {\bf 34} (1990),
  23-42.



\bibitem{BO}  A. Bondal and D. Orlov, {\em  Semiorthogonal decomposition for algebraic varieties}, preprint:math.AG/9506012.

\bibitem{CMZ} L. Costa and R.M. Mir\'o-Roig, {\em Tilting sheaves
on toric varieties}, Math. Z. {\bf 248} (2004), 849-865.

\bibitem{CMPAMS} L. Costa and R.M. Mir\'o-Roig, {\em Derived categories of projective
bundles},  Proc. Amer. Math. Soc.  {\bf 133}  (2005), 2533-2537.

\bibitem{CM} L. Costa and R.M. Mir\'o-Roig, {\em Derived category of Toric varieties
with small Picard number}, Preprint 2008.

\bibitem{CM1} L. Costa and R.M. Mir\'o-Roig, {\em Frobenius splitting and Derived category of Toric varieties}, To appear Illinois Journal of Mathematics.

\bibitem{bh} L. Borisov, Z. Hua, {\em On the conjecure of King for
smmoth toric Deligne-Mumford stacks}, math.AG/0801.2812v3.


\bibitem{Ewa} G. Ewald, {\em Combinatorial convexity and algebraic geometry}, Graduate Texts in Mathematics, {\bf 168} (1996) Springer-Verlag, New York.

\bibitem{Fu} W. Fulton, {\em Introduction to toric varieties},
Ann. of Math. Studies, Princeton, {\bf 131} (1993).


\bibitem{Ka} A.N. Rudakov et al. {\em Helices and vector bundles: Seminaire Rudakov}
Lecture Note Series, {\bf 148} (1990) Cambridge University Press.


\bibitem{HP} L. Hille, M. Perling, {\em A counterexample to King's
conjecture}, Compos. Math. {\bf 142} (2006),  1507-1521.


\bibitem{Kaw} Y. Kawamata, {\em Derived categories of toric varieties}, Michigan
Math. J. {\bf 54} (2006), 517-535.

\bibitem{Ko} M. Kontsevich, {\em Homological algebra of mirror symmetry},  Proceedings of the International Congress of Mathematicians, Vol. 1, 2 (1994),  120--139, Birkh{\"a}user, Basel, 1995.

\bibitem{MM} M. Lason, M. Michalek, {\em On the full, strongly exceptional collections on toric varieties with Picard number three}, arXiv:1003.2047, 2010.

\bibitem{Mateusz} M. Michalek, {\em Family of counterexamples to King's conjecture}, Preprint 2010.



\bibitem{Oda} T. Oda, {\em Convex Bodies and Algebraic Geometry}, Springer-Verlag
(1988).

\bibitem{Or} D.O. Orlov, {\em Projective bundles, monoidal transformations, and derived
categories of coherent sheaves}, Math. USSR Izv. {\bf 38} (1993),
133-141.


\end{thebibliography}
\end{document}